\documentclass[a4paper,11pt,twoside,reqno]{amsart}

\usepackage[utf8]{inputenc}
\usepackage[plainpages=false,pdfpagelabels=true]{hyperref}
\usepackage{amssymb,amsthm}
\usepackage[margin=1in]{geometry}

\newtheorem{Satz}{Theorem}[section]
\newtheorem{Prop}[Satz]{Proposition}
\newtheorem{Lem}[Satz]{Lemma}
\newtheorem{Cor}[Satz]{Corollary}
\newcommand{\vol}{{\operatorname{Area}}}
\theoremstyle{definition}

\newtheorem{Bem}[Satz]{Remark}
\allowdisplaybreaks[1]

\newcommand{\C}{\ensuremath{\mathbb{C}}}

\numberwithin{equation}{section}

\title{A note on twisted Dirac operators on closed surfaces}
\author{Volker Branding}
\date{\today}
\address{University of Vienna, Faculty of Mathematics\\
Oskar-Morgenstern-Platz 1, 1090 Vienna, Austria\\}
\email{volker.branding@univie.ac.at}
\begin{document}

\begin{abstract}
We derive an inequality that relates nodal set and eigenvalues of a class of twisted Dirac operators on closed surfaces
and point out how this inequality naturally arises as an eigenvalue estimate for the $\rm Spin^c$ Dirac operator.
This allows us to obtain eigenvalue estimates for the twisted Dirac operator appearing in the context of Dirac-harmonic maps and their extensions,
from which we also obtain several Liouville type results.
\end{abstract} 

\maketitle

\section{Introduction and results}
Throughout this note we assume that \((M,g)\) is a closed Riemannian spin surface with fixed spin structure.
On the spinor bundle \(\Sigma M\) we have a metric connection \(\nabla^{\Sigma M}\) induced from the Levi-Cevita connection 
and we fix a hermitian scalar product denoted by \(\langle\cdot,\cdot\rangle_{\Sigma M}\). 
Sections of the spinor bundle are called \emph{spinors}.
In addition, we have the Clifford multiplication of spinors with tangent vectors, denoted by \(X\cdot\psi\) for
\(X\in \Gamma(TM)\) and \(\psi\in\Gamma(\Sigma M)\). Clifford multiplication is skew-symmetric
\[
\langle X\cdot\psi,\xi\rangle_{\Sigma M}=-\langle \psi,X\cdot\xi\rangle_{\Sigma M}
\]
and satisfies the Clifford relation
\[
X\cdot Y\cdot\psi+Y\cdot X\cdot\psi=-2g(X,Y)\psi
\]
for \(X,Y\in \Gamma(TM)\) and \(\psi,\xi\in\Gamma(\Sigma M)\).

The Dirac operator is a first order differential operator acting on spinors and is given by
\[
D:=e_1\cdot\nabla^{\Sigma M}_{e_1}+e_2\cdot\nabla^{\Sigma M}_{e_2},
\]
where \(\{e_1,e_2\}\) is an orthonormal basis of \(TM\).
It is an elliptic differential operator, which is self-adjoint with respect to the \(L^2\)-norm.
On a compact Riemannian manifold the spectrum of the Dirac operator is discrete and
consists of both positive and negative eigenvalues. 
In general, there are only a few manifolds on which the spectrum of
the Dirac operator can be computed explicitly, see the book \cite{MR2509837}
for a current overview on the spectrum of the Dirac operator.

However, it is possible to estimate the spectrum. 
For instance, on a manifold of dimension \(n\geq 3\),
Friedrich \cite{MR600828} proved
that all eigenvalues \(\lambda\) of the Dirac operator satisfy
\[
\lambda^2\geq\frac{n}{4(n-1)}\inf_M R,
\]
where \(R\) denotes the scalar curvature of the manifold.
On closed surfaces Bär \cite{MR1162671} proved that
\[
\lambda^2\geq\frac{2\pi\chi(M)}{\vol(M,g)},
\]
where \(\chi(M)\) is the Euler characteristic of \(M\) and \(\vol(M,g)\)
denotes the area of \(M\).

Since \(\dim M=2\) and using the complex-volume form \(\omega_\C=ie_1\cdot e_2\),
we can decompose the spinor bundle into its positive and negative parts, that is
\(\Sigma M=\Sigma^+M\oplus\Sigma^-M\), where
\[
\Sigma^\pm M=\frac{1}{2}(1\pm\omega_\C)\cdot\Sigma M.
\]

The Dirac type operators that appear in theoretical physics usually do not act on sections of \(\Sigma M\)
but on sections of \(\Sigma M\otimes E\), where \(E\) is a given vector bundle. We call \(\Sigma M\otimes E\) 
a \emph{twisted spinor bundle}. The connections on \(\Sigma M\) and \(E\) induce a connection on \(\Sigma M\otimes E\)
and the same holds true for the scalar product on \(\Sigma M\otimes E\). Throughout this note
we will always take the real part of the scalar product on \(\Sigma M\otimes E\) turning it into
a Euclidean scalar product. It will be denoted by \(\langle\cdot,\cdot\rangle=\text{Re}(\langle\cdot,\cdot\rangle_{\Sigma M\otimes E})\).

The \emph{twisted Dirac operator} \(D^E\) maps sections of \(\Sigma M\otimes E\) to 
sections of \(\Sigma M\otimes E\) and is given by
\[
D^E:=e_1\cdot\nabla^{\Sigma M\otimes E}_{e_1}+e_2\cdot\nabla^{\Sigma M\otimes E}_{e_2}.
\]
Note that the Clifford multiplication only acts on the first factor of \(\Sigma M\otimes E\). 
The principal symbol of \(D^E\) can be computed as
\[
\sigma_1(D^E)=\sigma_1(D,\eta)\otimes id_E,\qquad \eta\in \Gamma(T^\ast M),
\]
hence this operator is still elliptic.
Moreover, \(D^E\) is formally self-adjoint whenever we have a metric connection on \(E\).
The square of the twisted Dirac operator satisfies the following Weitzenböck formula \cite[p.164]{MR1031992}
\begin{equation}
\label{weitzenboeck}
(D^E)^2=\nabla^\ast\nabla^{\Sigma M\otimes E}+\frac{R}{4}+\frac{1}{2}\sum_{i,j=1}^2e_i\cdot e_j\cdot R^E(e_i,e_j),
\end{equation}
where \(R^E\) denotes the curvature endomorphism of
the twist bundle \(E\).

For more background material on spin geometry and the Dirac operator we refer the reader to the books \cite{MR1476425} and \cite{MR1031992}.

We will always assume that we have a metric connection on \(E\).
A twisted spinor \(\psi\in\Gamma(\Sigma M\otimes E)\) is called \emph{eigenspinor} of the twisted Dirac operator \(D^E\) with eigenvalue \(\lambda\) if it satisfies
\[
D^E\psi=\lambda\psi.
\]
The twisted Dirac operator \(D^E\) interchanges positive and negative twisted spinors, more precisely
\(D^E:\Sigma^\pm M\otimes E\to\Sigma^\mp M\otimes E\), whereas its square \((D^E)^2:\Sigma^\pm M\otimes E\to\Sigma^\pm M\otimes E\)
preserves the grading.

Of course, one cannot expect to compute the spectrum of the twisted Dirac operator explicitly.
Nevertheless some eigenvalue estimates exist, see for example \cite{MR2737611,MR2536121,MR2521901,MR2241738,MR2231389}.

In this note we will analyze the nodal set of spinors
\begin{align*}
\mathcal{N}(\psi):=\{p\in M\mid\psi(p)=0\}
\end{align*}
that solve semilinear elliptic equations of first order.
By the main result of \cite{MR1714341} we know that this set is discrete.

In this note we generalize an inequality relating zero-set
and eigenvalues of the classical Dirac operator obtained in \cite{Branding2017}, see also \cite{MR3210753},
to a large class of twisted spinors.
Although this generalization is rather straightforward we will provide several applications
of the main result to twisted Dirac operators that arise in quantum field theory.

\begin{Satz}
\label{main-result}
Let \((M,g)\) be a closed spin surface with fixed spin structure and let \(\Sigma M\otimes E\)
be the spinor bundle twisted by a complex line bundle \(E\) with metric connection.
Then all eigenvalues of \(D^E\) satisfy
\begin{equation}
\label{eigenvalue-estimate}
\lambda^2\geq\frac{2\pi\chi(M)}{\vol(M,g)}-\frac{4\pi}{\vol(M,g)}\int_M|c_1(E)|dM+\frac{4\pi N(\psi)}{\vol(M,g)},
\end{equation}
where \(\chi(M)\) is the Euler characteristic of \(M\) and \(c_1(E)\) the first Chern-form of the complex line bundle \(E\).
Moreover, \(N(\psi)\) denotes the sum of the orders of the zeros of the eigenspinor \(\psi\)
\begin{equation}
\label{nodal-set}
N(\psi)=\sum_{p\in M,|\psi|(p)=0}n_p,
\end{equation}
where \(n_p\) represents the order of the corresponding zero.
\end{Satz}

In the case of \(\lambda=0\) we can explicitly calculate the nodal set of \(\psi\) similar to
\cite[Proposition 3.2]{Branding2017}.

\begin{Satz}
\label{main-result-2}
Let \((M,g)\) be a closed spin surface with fixed spin structure and let \(\Sigma M\otimes E\)
be the spinor bundle twisted by a complex line bundle \(E\) with metric connection.
Let \(\psi\in\Gamma(\Sigma^+M\otimes E)\) be a solution of \(D^E\psi=0\). Then the following equality holds
\begin{equation}
\label{nodal-set-harmonic}
N_0(\psi)=-\frac{\chi(M)}{2}-\deg(E),
\end{equation}
where \(\chi(M)\) is the Euler characteristic of \(M\), \(\deg(E)\) the degree of the complex line bundle \(E\)
and \(N_0(\psi)\) denotes the sum of the orders of the zeros of \(\psi\).
\end{Satz}

\begin{Bem}
If \(\psi\in\Gamma(\Sigma^-M\otimes E)\) is a solution of \(D^E\psi=0\), then 
the nodal set satisfies \(N_0(\psi)=-\frac{\chi(M)}{2}+\deg(E)\).
\end{Bem}

\begin{Cor}
Of course, \eqref{eigenvalue-estimate} can also be interpreted as an estimate on the nodal set of an eigenspinor \(\psi\), that is
\[
N(\psi)\leq\frac{\vol(M,g)\lambda^2}{4\pi}-\frac{\chi(M)}{2}+\int_M|c_1(E)|dM.
\]
\end{Cor}

The two main results of this article can also be interpreted from a more geometric perspective, 
namely as an eigenvalue estimate for the $\rm Spin^c$ Dirac operator.
\footnote{The author would like to thank Georges Habib and Nicolas Ginoux for pointing out the relation to $\rm Spin^c$ structures.}

A Riemannian manifold admits a $\rm Spin^c$ structure if and only if there exists a complex
line bundle \(E\) such that \([c_1(E)]_{mod~2}=\omega_2(M)\), where \(\omega_2(M)\) is the second
Stiefel-Whintey class of the manifold \(M\).
In particular, every Riemannian spin surface also has a $\rm Spin^c$ structure.
We can think of the Dirac operator associated to a $\rm Spin^c$ structure as the twisted Dirac operator
acting on spinors which are twisted by a line bundle.
For more details on $\rm Spin^c$ structures see \cite[Appendix D]{MR1031992} and \cite{MR2959436,MR2661160}
for eigenvalue estimates of the  $\rm Spin^c$ Dirac operator.

For this reason the eigenvalue estimate \eqref{eigenvalue-estimate} can also be interpreted as
an eigenvalue estimate of the Dirac operator associated to any $\rm Spin^c$ structure
and \eqref{nodal-set-harmonic} as an estimate on the nodal set of a spinor that is
in the kernel of the $\rm Spin^c$ Dirac operator.

Theorem \ref{main-result} is similar to Theorem 3.4 in \cite{MR2959436}
taking also into account the zero set of \(\psi\). However, note that
the constant in front of the second term in \eqref{eigenvalue-estimate} is different
since we are estimating the Chern-form instead of the curvature two-form of the line bundle 
associated to the $\rm Spin^c$ structure.

In addition, Theorem 3.4 in \cite{MR2959436} also discusses the equality case.
It turns out that one possible case of equality in \eqref{eigenvalue-estimate} is achieved 
under the same assumptions. More precisely, if \(\psi\) is 
a $\rm Spin^c$-Killing spinor (which does not have any zeros), that is a solution of
\begin{align*}
\label{killing-spinorc}
\nabla^{\Sigma M\otimes E}_X\psi=-\frac{\lambda}{2}X\cdot\psi,\qquad X\in\Gamma(TM) 
\end{align*}
and satisfies \(\omega_\C\cdot\psi=-\frac{|c_1(E)|}{c_1(E)}\psi\) at the same time
then we obtain equality in \eqref{eigenvalue-estimate}.

This article is organized as follows: In the second section we will give the proofs of the two main theorems,
the third section provides several applications of the main results to twisted Dirac operators
originating in quantum field theory.

\section{Proof of the main Theorems}
By the main result of \cite{MR1714341} we know that on a two-dimensional manifold the zero-set of
twisted eigenspinors is discrete. In the following we will make use of the energy-momentum tensor \(T^E(X,Y)\),
which is given by
\[
T^E(X,Y):=\langle X\cdot\nabla^{\Sigma M\otimes E}_Y\psi+Y\cdot\nabla^{\Sigma M\otimes E}_X\psi,\psi\rangle.
\]
This tensor arises if one varies the functional \(E(\psi)=\int_M\langle\psi,D^E\psi\rangle dM\) with respect to the metric.
The following Lemma is a generalization of Lemma 5.1 from \cite{MR1738150}.
Note that it is essential that we take the real part of the scalar product on \(\Sigma M\otimes E\) in order for the next Lemma to hold. 
\begin{Lem}
\label{lem-modified-connection}
Suppose that \(\psi\in\Gamma(\Sigma M\otimes E)\) is nowhere vanishing. Then the following inequality holds
\begin{equation}
\label{inequality-modified-connection}
\frac{\langle\psi,(D^E)^2\psi\rangle}{|\psi|^2}\geq\frac{R}{4}+\frac{1}{2}\sum_{i,j=1}^2\frac{\langle e_i\cdot e_j\cdot R^E(e_i,e_j)\psi,\psi\rangle}{|\psi|^2}+\frac{|T^E|^2}{4|\psi|^4}-\Delta\log|\psi|
-\frac{\langle D^E\psi,d(\log|\psi|^2)\cdot\psi\rangle}{|\psi|^2}.
\end{equation}
\end{Lem}
\begin{proof}
We set
\[
\tilde{\nabla}^{\Sigma M\otimes E}_X\psi:=\nabla_X^{\Sigma M\otimes E}\psi-2\alpha(X)\psi-\beta(X)\cdot\psi-X\cdot\alpha\cdot\psi
\]
with a one-form \(\alpha\) and a symmetric \((1,1)\)-tensor \(\beta\) given by
\[
\alpha:=\frac{d|\psi|^2}{2|\psi|^2}, \qquad \beta:=-\frac{T^E(\cdot,\cdot)}{2|\psi|^2}.
\]
By a direct computation we then find summing over repeated indices
\begin{align*}
|\tilde{\nabla}^{\Sigma M\otimes E}\psi|^2=&|\nabla^{\Sigma M\otimes E}\psi|^2+2|\alpha|^2|\psi|^2+|\beta|^2|\psi|^2
-4\alpha({e_i})\langle\nabla^{\Sigma M\otimes E}_{e_i}\psi,\psi\rangle \\
&+2\langle\beta(e_i)\cdot\nabla^{\Sigma M\otimes E}_{e_i}\psi,\psi\rangle
+2\langle D^E\psi,\alpha\cdot\psi\rangle.
\end{align*}
Moreover, using \eqref{weitzenboeck}, we obtain
\begin{align*}
|\nabla^{\Sigma M\otimes E}\psi|^2=&\langle\psi,(D^E)^2\psi\rangle-\frac{R}{4}|\psi|^2-\frac{1}{2}\sum_{i,j=1}^2\langle e_i\cdot e_j\cdot R^E(e_i,e_j)\psi,\psi\rangle+\frac{1}{2}\Delta|\psi|^2, \\
\alpha({e_i})\langle\nabla^{\Sigma M\otimes E}_{e_i}\psi,\psi\rangle=&|\alpha|^2|\psi|^2=\frac{|d|\psi|^2|^2}{4|\psi|^2}, \\
\langle\beta(e_i)\cdot\nabla^{\Sigma M\otimes E}_{e_i}\psi,\psi\rangle=&-\frac{|T^E|^2}{4|\psi|^2}.
\end{align*}
Thus, we get after summing over repeated indices
\begin{align*}
0\leq|\tilde{\nabla}^{\Sigma M\otimes E}\psi|^2=&\langle\psi,(D^E)^2\psi\rangle-\frac{R}{4}|\psi|^2
-\frac{1}{2}\langle e_i\cdot e_j\cdot R^E(e_i,e_j)\psi,\psi\rangle
+\frac{1}{2}\Delta|\psi|^2-\frac{|d|\psi|^2|^2}{2|\psi|^2}\\
&-\frac{|T^E|^2}{4|\psi|^2} 
+2\langle D^E\psi,\alpha\cdot\psi\rangle,
\end{align*}
which yields the result.
\end{proof}

Throughout the article we will make use of the following fact:
If the zero set of \(|\psi|\) is discrete and \(|\psi|\) does not vanish identically, 
then the following equality holds
\begin{equation}
\label{laplace-log}
\int_M\Delta\log|\psi|dM=-2\pi\sum_{p\in M,|\psi|(p)=0}n_p,
\end{equation}
where \(n_p\) is the order of \(|\psi|\) at \(p\).
A proof can be found in \cite{MR1474501}, see also \cite[Lemma 2.2]{Branding2017}.

\begin{proof}[Proof of Theorem \ref{main-result}]
To complete the proof of Theorem \ref{main-result} we analyze the term arising from the curvature of the twist bundle \(E\) in Lemma \ref{lem-modified-connection}.
Using the skew symmetry of the curvature endomorphism and the Clifford multiplication we get
\[
\sum_{i,j=1}^2e_i\cdot e_j \cdot R^E(e_i,e_j)\psi=2e_1\cdot e_2\cdot R^E(e_1,e_2)\psi
=-2i\omega_\C\cdot R^E(e_1,e_2)\psi=4\pi c_1(E)\omega_\C\cdot\psi,
\]
where \(c_1(E)\) denotes the first Chern-form of the complex line bundle \(E\), see \cite{MR0440554}, p.303 ff. for more details.
Hence, we find
\[
\frac{1}{2}\sum_{i,j=1}^2\frac{\langle e_i\cdot e_j\cdot R^E(e_i,e_j)\psi,\psi\rangle}{|\psi|^2}=2\pi c_1(E)\frac{\langle\omega_\C\cdot\psi,\psi\rangle}{|\psi|^2}.
\]
Now, we apply Lemma \ref{lem-modified-connection} in the case that \(\psi\) is an eigenspinor of \(D^E\).
We can estimate the energy momentum tensor by \(|T|^2\geq 2\lambda^2|\psi|^4\).
Moreover, the last term on the right hand side of \eqref{inequality-modified-connection} vanishes
since it is both purely real and purely imaginary.
Thus, from \eqref{inequality-modified-connection} we obtain
\begin{align*}
\lambda^2\geq& K+4\pi c_1(E)\frac{\langle\omega_\C\cdot\psi,\psi\rangle}{|\psi|^2}-2\Delta\log|\psi|
\\ \geq& K-4\pi |c_1(E)|-2\Delta\log|\psi|, 
\end{align*}
where \(K\) denotes the Gaussian curvature of \(M\). 
By integrating over the surface \(M\) and using \eqref{laplace-log} we obtain Theorem \ref{main-result}.
\end{proof}

In order to prove Theorem \ref{main-result-2} we will make use of the following Bochner-formula.

\begin{Lem}
Let \(\psi\in\Gamma(\Sigma^+M\otimes E)\) be a solution of \(D^E\psi=0\).
Then the following equality holds
\begin{align}
\Delta\log|\psi|=\frac{K}{2}+2\pi c_1(E).
\end{align}
\end{Lem}
\begin{proof}
By assumption we have \(\psi\in\Gamma(\Sigma M^+\otimes E)\), thus we find
\[
\frac{\langle\omega_\C\cdot\psi,\psi\rangle}{|\psi|^2}=1.
\]
By a direct calculation using the Weitzenböck formula \eqref{weitzenboeck} we then get
\begin{align}
\label{identity-d}
\Delta\log|\psi|=\frac{K}{2}+\frac{|\nabla^{\Sigma M\otimes E}\psi|^2}{|\psi|^2}-\frac{1}{2}\frac{\big|d|\psi|^2\big|^2}{|\psi|^4}+2\pi c_1(E).
\end{align}
Now, we use the same arguments as in the proof of \cite[Proposition 3.2]{Branding2017}.
Since \(\Sigma^+M\otimes E\) is a complex line bundle, we can write \(\nabla^{\Sigma M\otimes E}_{e_j}\psi=f_j\psi\)
for some complex-valued function \(f_j\) away from its zero-set.
By a direct calculation we obtain the following identities
\begin{align}
\label{identity-a}\big|d|\psi|^2\big|^2=&4\sum_{j=1}^2|\operatorname{Re}f_j|^2|\psi|^4, \\
\label{identity-b}|\nabla^{\Sigma M\otimes E}\psi|^2=&\sum_{j=1}^2|f_j|^2|\psi|^2, \\
\label{identity-c}|D^E\psi|^2=&\big((\operatorname{Re}f_1+\operatorname{Im}f_2)^2+(\operatorname{Re}f_2+\operatorname{Im}f_1)^2\big)|\psi|^2.
\end{align}
By assumption \(D^E\psi=0\), consequently we get from \eqref{identity-c} that \(\operatorname{Re}f_1=-\operatorname{Im}f_2\)
and \(\operatorname{Re}f_2=-\operatorname{Im}f_1\). Inserting this into both \eqref{identity-a} and \eqref{identity-b}
we obtain
\begin{align*}
|\nabla^{\Sigma M\otimes E}\psi|^2=\frac{1}{2}\frac{\big|d|\psi|^2\big|^2}{|\psi|^2}.
\end{align*}
Combining this identity with \eqref{identity-d} completes the proof.
\end{proof}

\begin{proof}[Proof of Theorem \ref{main-result-2}]
  By integrating the identity
\begin{align*}
\Delta\log|\psi|=\frac{K}{2}+2\pi c_1(E)
\end{align*}
over \(M\) and recalling that
the degree of a vector bundle is given by
\[
\deg(E)=\int_Mc_1(E)dM,
\]
we obtain the desired result. 
\end{proof}

\section{Applications of the main results}
\subsection{Dirac-harmonic maps from surfaces}
In this section we focus on Dirac-harmonic maps \cite{MR2262709} from surfaces.
Thus, let \(N\) be another Riemannian manifold, let \(\phi\colon M\to N\)
be a map and \(E\) be the pull-back bundle \(\phi^\ast TN\) over \(M\).
We are considering so-called \emph{vector spinors}, that is \(\psi\in\Gamma(\Sigma M\otimes \phi^\ast TN)\).
Dirac-harmonic maps are critical points of the energy functional (with the first term being the usual Dirichlet energy)
\[
E(\phi,\psi)=\frac{1}{2}\int_M(|d\phi|^2+\langle\psi,D^{\phi^\ast TN}\psi\rangle)dM
\]
and satisfy the Euler-Lagrange equations
\[
\tau(\phi)=\frac{1}{2}\sum_{i=1}^2R^N(\psi,e_i\cdot\psi)d\phi(e_i),\qquad D^{\phi^\ast TN}\psi=0.
\]
Here, \(\tau(\phi)\) is the tension field of the map \(\phi\) and \(R^N\) the curvature tensor on \(N\).
Since the vector spinors \(\psi\) behave like tangent vectors on \(N\) they can be inserted into the curvature tensor \(R^N\).
This can be understood as follows: Locally we can expand a vector spinor \(\psi\in\Gamma(\Sigma M\otimes \phi^\ast TN)\) as
\begin{align*}
\psi=\sum_{\alpha=1}^{\dim N}\psi^\alpha\otimes\frac{\partial}{\partial y^\alpha},
\end{align*}
where \(\psi^\alpha,\alpha=1,\ldots\dim N\) are sections in the spinor bundle \(\Sigma M\)
and \(\frac{\partial}{\partial y^\alpha},\alpha=1,\ldots\dim N\) is a local basis of \(\phi^\ast TN\).
In terms of local coordinates \(x\) on \(M\) and \(y\) on \(N\) we can write
\[
\sum_{i=1}^2R^N(e_i\cdot\psi,\psi)d\phi(e_i):=
\sum_{i=1}^2\sum_{\alpha,\beta,\gamma,\delta=1}^{\dim N} R^\alpha_{~\beta\gamma\delta}\frac{\partial}{\partial y^\alpha}\langle\psi^\gamma,e_i\cdot\psi^\delta\rangle_{\Sigma M}\frac{\partial\phi^\beta}{\partial x_i},
\]
where we used Latin indices on \(M\) and Greek indices on \(N\). Here \(R^\alpha_{~\beta\gamma\delta}\) are the components of
the curvature tensor on the target \(N\).

\begin{Prop}
\label{prop-dhmaps}
Let \(M\) be a closed Riemannian spin surface with fixed spin structure and \(N\) another Riemannian manifold.
Suppose that \((\phi,\psi)\) is a smooth non-trivial Dirac-harmonic map from \(M\) to \(N\).
Then the following inequality holds
\begin{equation}
\label{inequality-dhmap}
|R^N|_{L^\infty}\int_M|d\phi|^2dM\geq\pi\chi(M)+2\pi N_0(\psi),
\end{equation}
where \(N_0(\psi)\) denotes the sum of the order of the zeros of \(\psi\).
\end{Prop}

\begin{proof}
First of all, we note that 
\[
\sum_{i,j=1}^2\frac{1}{2}\frac{\langle e_i\cdot e_j\cdot R^N(d\phi(e_i),d\phi(e_j))\psi,\psi\rangle}{|\psi|^2}\geq-|R^N|_{L^\infty}|d\phi|^2.
\]
Using \eqref{inequality-modified-connection} with \(E=\phi^\ast TN\) and the assumption \(D^{\phi^\ast TN}\psi=0\) we get
\[
|R^N|_{L^\infty}|d\phi|^2\geq\frac{R}{4}-\Delta\log|\psi|.
\]
The result then follows by integration.
\end{proof}
Making use of Proposition \ref{prop-dhmaps} we obtain the following:

\begin{Cor}
Let \(M\) be a closed Riemannian spin surface and \(N\) another Riemannian manifold.
Suppose that \((\phi,\psi)\) is a smooth non-trivial Dirac-harmonic map from \(M\) to \(N\).
If \(\chi(M)>0\) and the energy of \(\phi\) is sufficiently small, that is
\[
\int_M|d\phi|^2dM<\epsilon
\]
for some small \(\epsilon>0\) then the pair \((\phi,\psi)\) is trivial.
\end{Cor}
\begin{proof}
The triviality of \(\psi\) follows from \eqref{inequality-dhmap}, for the vanishing of the map \(\phi\) 
see \cite[Lemma 4.9]{MR3333092}.
\end{proof}

If both \(M=N=S^2\) it was proven in \cite{MR2496649} that the Euler-Lagrange equations for Dirac-harmonic maps
decouple, that is
\[
\tau(\phi)=0=\frac{1}{2}\sum_{i=1}^2R^N(\psi,e_i\cdot\psi)d\phi(e_i),\qquad D^{\phi^\ast TN}\psi=0.
\]
Making use of this fact, all non-trivial Dirac-harmonic maps between spheres could be classified.
Moreover, it is well-known that harmonic maps between two-spheres are holomorphic or anti-holomorphic maps with energy
\[
\frac{1}{2}\int_M|d\phi|^2dM=E(\phi)=4\pi|\deg(\phi)|.
\]
Hence, we obtain the following inequality for uncoupled Dirac-harmonic maps between spheres
\begin{equation*}
4|R^N|_{L^\infty}|\deg(\phi)|\geq 1+N_0(\psi).
\end{equation*}

\subsection{Dirac-harmonic maps between surfaces}
In this section we apply Theorem \ref{main-result-2} in the case of Dirac-harmonic maps between surfaces.
More precisely, we assume that \(N\) is a closed, connected and oriented Riemannian surface.
Thus, we can think of \(TN\) as a complex line bundle.
The following equality was already proven in \cite{MR2496649}
giving rise to a structure theorem for Dirac-harmonic maps between closed surfaces.

\begin{Prop}
Let \(E=\phi^\ast TN\) and suppose that \(\psi\in\Gamma(\Sigma^+ M\otimes\phi^\ast TN)\)
is a solution of \(D^{\phi^\ast TN}\psi=0\).
Then the following formula holds
\begin{equation}
\label{estimate-dirac-harmonic}
N_0(\psi)=-\frac{\chi(M)}{2}-\deg(\phi)\chi(N),
\end{equation}
where \(\deg(\phi)\) denotes the degree of the map \(\phi\) and \(\chi(N)\) is the Euler-characteristic
of \(N\). 
\end{Prop}

\begin{proof}
We apply Theorem \ref{main-result-2} in the case that \(E=\phi^\ast TN\).
By the properties of the first Chern-form we find
\[
\deg(\phi^\ast TN)=\int_Mc_1(\phi^\ast TN)dM=\int_M\phi^\ast c_1(TN)dM=\frac{1}{2\pi}\int_M\phi^\ast K^NdM.
\]
Using the degree theorem we get
\[
\int_M\phi^\ast K^NdM=\deg(\phi)\int_NK^NdN=2\pi\deg(\phi)\chi(N),
\]
which proves the result.
\end{proof}

Equality \eqref{estimate-dirac-harmonic} still holds if we consider \emph{Dirac-harmonic maps with torsion} between surfaces,
which were introduced in \cite{MR3493217}. In this case one considers a metric connection
with torsion on \(TN\). On a surface a metric connection with torsion \(\nabla^{\scriptscriptstyle Tor}\) is given by
\begin{align*}
\nabla_X^{\scriptscriptstyle Tor}Y=\nabla_X^{\scriptscriptstyle LC}Y+\langle X,Y\rangle V-\langle V,Y\rangle X,
\end{align*}
where \(X,Y,V\) are vector fields on \(N\) and \(\nabla^{\scriptscriptstyle LC}\) denotes the Levi-Cevita connection.

Moreover, for a metric connection with torsion on a surface the Gauss curvatures satisfy
\begin{align}
\label{curvature-torsion}
K^{\scriptscriptstyle Tor}=K^{\scriptscriptstyle LC}+\operatorname{div}^{LC} V,
\end{align}
where \(\operatorname{div}^{LC}\) denotes the divergence with respect to the Levi-Cevita connection.

\begin{Prop}
Let \(E=\phi^\ast TN\) equipped with a metric connection with torsion 
and suppose that \(\psi\in\Gamma(\Sigma^+ M\otimes\phi^\ast TN)\) is a solution of \(D^{\phi^\ast TN}\psi=0\).
Then the following formula holds
\begin{equation}
\label{estimate-dirac-harmonic-torsion}
N_0(\psi)=-\frac{\chi(M)}{2}-\deg(\phi)\chi(N),
\end{equation}
where \(\deg(\phi)\) denotes the degree of the map \(\phi\) and \(\chi(N)\) is the Euler-characteristic
of \(N\). 
\end{Prop}
\begin{proof}
The proof is essentially the same as the previous one:
We apply Theorem \ref{main-result-2} and make use of \eqref{curvature-torsion}, the additional
divergence term arising from the torsion vanishes after integration.

The statement is also clear from a more general point of view: 
The degree of the vector bundle \(\phi^\ast TN\) does not depend on the chosen
connection and thus the torsion does not contribute.
\end{proof}

\subsection{Spinor-valued one-forms}
For \(E=T^\ast M\) we call sections in \(\Sigma M\otimes T^\ast M\)
spinor-valued one-forms. These appear in the context of the Rarita-Schwinger operator,
however one needs to introduce a certain projection in order to relate the twisted Dirac operator
on \(\Sigma M\otimes T^\ast M\) to the Rarita-Schwinger operator,
see \cite{MR3289848}, section 2.6, \cite{MR1891206} and \cite{MR1129331} for more details.
Nevertheless, we apply our main result in the case that \(E=T^\ast M\).
\begin{Prop}
Let \(E=T^\ast M\) and suppose that \(\psi\in\Gamma(\Sigma^+ M\otimes T^\ast M)\)
is a solution of \(D^{T^\ast M}\psi=0\).
Then the following formula holds
\begin{equation}
N_0(\psi)=\frac{1}{2}\chi(M).
\end{equation}
\end{Prop}
\begin{proof}
This follows directly from Theorem \ref{main-result-2} by using
\(
\deg(T^\ast M)=-\deg(TM)=-\chi(M).
\)
\end{proof}

\subsection{Dirac-harmonic maps with curvature term}
\emph{Dirac-harmonic maps with curvature term} arise as critical points of the energy functional for Dirac-harmonic maps 
taking into account an additional curvature term, see \cite{MR3333092, MR2370260} for more details.
We again consider \(E=\phi^\ast TN\) and the manifold \(N\) may have arbitrary dimension.
Dirac-harmonic maps with curvature term arise as critical points of the energy functional
\begin{align*}
E_c(\phi,\psi)=\frac{1}{2}\int_M(|d\phi|^2+\langle\psi,D^{\phi^\ast TN}\psi\rangle-\frac{1}{6}\langle R^N(\psi,\psi)\psi,\psi\rangle)dM,
\end{align*}
where the spinors are contracted as
\[
\langle R^N(\psi,\psi)\psi,\psi\rangle_{\Sigma M\otimes\phi^\ast TN}=\sum_{\alpha,\beta,\gamma,\delta=1}^{\dim N}R_{\alpha\beta\gamma\delta}
\langle\psi^\alpha,\psi^\gamma\rangle_{\Sigma M}\langle\psi^\beta,\psi^\delta\rangle_{\Sigma M}.
\]
Here, \(R_{\alpha\beta\gamma\delta}\) are the components of the Riemann curvature tensor on \(N\).
Dirac-harmonic maps with curvature term are solutions of the coupled system
\begin{align}
\nonumber\tau(\phi)=&\frac{1}{2}\sum_{i=1}^2R^N(\psi,e_i\cdot\psi)d\phi(e_i)-\frac{1}{12}\langle(\nabla R^N)^\sharp(\psi,\psi)\psi,\psi\rangle, \\
\label{psi-dhmap-curvature}D^{\phi^\ast TN}\psi=&\frac{1}{3}R^N(\psi,\psi)\psi.
\end{align}
The curvature term on the right hand side of \eqref{psi-dhmap-curvature} is defined as follows:
In terms of local coordinates \(y\) on \(N\) we can write
\begin{align*}
R^N(\psi,\psi)\psi:=\sum_{\alpha,\beta,\gamma,\delta=1}^{\dim N} R^\alpha_{~\beta\gamma\delta}\frac{\partial}{\partial y^\alpha}\langle\psi^\beta,\psi^\delta\rangle_{\Sigma M}\psi^\gamma,
\end{align*}
where \(\psi^\alpha,\alpha=1,\ldots,\dim N\) are sections in \(\Sigma M\).

\begin{Prop}
Let \((\phi,\psi)\) be a smooth Dirac-harmonic map with curvature term
and energy \(E(\phi,\psi)=\int_M(|d\phi|^2+|\psi|^4)dM\). Then the following inequality holds
\begin{equation}
\label{inequality-dhmap-curvature-psi}
C_NE(\phi,\psi)\geq\pi\chi(M)+2\pi N(\psi)
\end{equation}
with the constant \(C_N:=\max\{\frac{|R^N|_{L^\infty}^2}{6},|R^N|_{L^\infty}\}\) and \(N(\psi)\) is defined as in \eqref{nodal-set}.
\end{Prop}
\begin{proof}
First of all we note, using that \(\psi\) is a solution of \eqref{psi-dhmap-curvature}, that
\begin{align*}
\langle\psi,(D^{\phi^\ast TN})^2\psi\rangle=&\frac{1}{3}\langle\psi,D^{\phi^\ast TN}(R^N(\psi,\psi)\psi)\rangle \\
=&\frac{1}{3}\langle\psi,(\nabla(R^N(\psi,\psi)))\cdot\psi)\rangle+\frac{1}{3}\langle\psi,R^N(\psi,\psi)D^{\phi^\ast TN}\psi\rangle \\
=&\frac{1}{9}|R^N(\psi,\psi)\psi|^2,
\end{align*}
where we used the skew-symmetry of the Clifford multiplication. Moreover, we can estimate the energy momentum tensor as follows
\[
|T^{\phi^\ast TN}|^2\geq 2|\langle\psi,D^{\phi^\ast TN}\psi\rangle|^2=\frac{2}{9}|\langle\psi,R^N(\psi,\psi)\psi\rangle|^2.
\]
Note that the last term on the right hand side of \eqref{inequality-modified-connection} vanishes
since it is both purely real and purely imaginary when \(\psi\) solves \eqref{psi-dhmap-curvature}.
Hence from \eqref{inequality-modified-connection} with \(E=\phi^\ast TN\) we get
\[
\frac{1}{9}\frac{|R^N(\psi,\psi)\psi|^2}{|\psi|^2}\geq\frac{K}{2}+\frac{1}{2}\sum_{i,j=1}^2\frac{\langle e_i\cdot e_j\cdot R^N(d\phi(e_i),d\phi(e_j)\psi,\psi\rangle}{|\psi|^2}
+\frac{|\langle\psi,R^N(\psi,\psi)\psi\rangle|^2}{18|\psi|^4}-\Delta\log|\psi|.
\]
Proceeding as in the proof of Proposition \ref{prop-dhmaps} we obtain
\[
C_N(|d\phi|^2+|\psi|^4)\geq\frac{K}{2}-\Delta\log|\psi|,
\]
where \(C_N:=\max\{\frac{|R^N|_{L^\infty}^2}{6},|R^N|_{L^\infty}\}\). 
The claim then follows by integration.
\end{proof}
The above Proposition allows us to draw the following conclusion:

\begin{Cor}
Suppose \((\phi,\psi)\) is a smooth Dirac-harmonic map with curvature term.
If \(\chi(M)>0\) and if the energy is sufficiently small, that is \(E(\phi,\psi)<\epsilon\) for some small \(\epsilon>0\),
then we get a contradiction from \eqref{inequality-dhmap-curvature-psi} forcing \(\psi\) to be trivial.
This has already been proven in \cite[Lemma 4.9]{MR3333092} with the help of the Sobolev embedding theorem.
However, in \eqref{inequality-dhmap-curvature-psi} all the constants are explicit.
\end{Cor}

\par\medskip
\textbf{Acknowledgements:}
The author gratefully acknowledges the support of the Austrian Science Fund (FWF): P 30749
through the project ``Geometric variational problems from string theory''.
\bibliographystyle{plain}
\bibliography{mybib}
\end{document}